\documentclass{article}
\usepackage{latexsym}
\newtheorem{llemma}{Lemma}[section]
\newtheorem{prop}[llemma]{Proposition}
\newtheorem{pf}[llemma]{Proof}
\newtheorem{exmp}[llemma]{Example}
\newtheorem{thm}[llemma]{Theorem}

\newtheorem{defn}[llemma]{Definition}
\newtheorem{rem}[llemma]{Remark}

\newtheorem{nnote}[llemma]{Note}

\begin{document}
%*********************************************
% TITLE
%\title{\textbf{A new approaches to soft sets with interval neutrosophic sets}}
\title{\textbf{Interval-valued neutrosophic soft sets and its decision making}}
\author{Irfan Deli\\ \\
%Muallim R{\i}fat Faculty of Education,\\
 Kilis 7 Aral{\i}k University, 79000 Kilis, Turkey, \\
              (irfandeli@kilis.edu.tr)}
\date{}

\maketitle
\begin{abstract}
In this paper, the notion of the interval valued neutrosophic soft
sets ($ivn-$soft sets) is defined which is a combination of an
interval valued neutrosophic sets \cite{wan-05}  and a soft sets
\cite{mol-99}. Our $ivn-$soft sets generalizes the concept of the
soft set, fuzzy soft set, interval valued fuzzy soft set,
intuitionistic fuzzy soft set, interval valued intuitionistic fuzzy
soft set and neutrosophic soft set. Then, we introduce some
definitions and operations on $ivn-$soft sets sets. Some properties
of $ivn-$soft sets which are connected to operations have been
established. Also, the aim of this paper is to investigate the
decision making based on $ivn-$soft sets by level soft sets.
Therefore, we develop a decision making methods and then give a
example to illustrate the developed approach.
\end{abstract}

\section{Introduction}
%---------------------------------------------------------------
Many fields deal with the uncertain data may not be successfully
modeled by the classical mathematics, since concept of uncertainty
is too complicate and not clearly defined object. But they can be
modeled a number of different approaches including the probability
theory, fuzzy set theory \cite{zad-65}, rough set theory
\cite{paw-82}, neutrosophic set theory \cite{sma-05} and some other
mathematical tools. This theory have been applied in many real
applications to handle uncertainty. In 1999, Molodtsov \cite{mol-99}
succesfully proposed a completely new theory so-called {\em soft set
theory} by using classical sets because its been pointed out that
soft sets are not appropriate to deal with uncertain and fuzzy
parameters. The theory is a relatively new mathematical model for
dealing with uncertainty from a parametrization point of view.

After Molodtsov, there has been a rapid growth of interest in soft
sets and their various applications such as; algebraic structures
(e.g.\cite{acar-10, akt-07, ayg-09, zhan-10}), ontology
(e.g.\cite{jia-11a}), optimization (e.g.\cite{kov-07}), lattice
(e.g.\cite{cag-12,nag-11, yýl-13}), topology
(e.g.\cite{cag-12d,min-11, sha-11}), perron integration, data
analysis and operations research (e.g.\cite{mol-99, mol-04}), game
theory (e.g.\cite{deli-13a, deli-13b, mol-99}), clustering
(e.g.\cite{awa-10, mam-13}), medical diagnosis (e.g.\cite{kha-10,
yük-13}), and decision making under uncertainty
(e.g.\cite{cag-09a,fen-08b, maj-02}). In recent years, many
interesting applications of soft set theory have been expanded by
embedding the ideas of fuzzy sets (e.g. \cite{cag-09c,
deli-12a,deli-12b,fen-08b,maj-01a}), rough sets (e.g. \cite{fen-10})
and intuitionistic fuzzy sets (e.g. \cite{jia-11,maj-01c}), interval
valued intuitionistic fuzzy (e.g. \cite{jia-10,zha-13}),
Neutrosophic (e.g. \cite{maj-13,maj-12}).

Intuitionistic fuzzy sets can only handle incomplete information
because the sum of degree true, indeterminacy and false is one in
intuitionistic fuzzy sets. But neutrosophic sets can handle the
indeterminate information and inconsistent information which exists
commonly in belief systems in neutrosophic set since indeterminacy
is quantified explicitly and truth-membership,
indeterminacy-membership and falsity-membership are independent. It
is mentioned in \cite{wan-05}. Therefore, Maji firstly proposed
neutrosophic soft sets with operations, which is free of the
difficulties mentioned above, in \cite{maj-13}. He also, applied to
decision making problems in \cite{maj-12}. After Maji, the studies
on the neutrosophic soft set theory have been studied increasingly
(e.g. \cite{bro-13,bro-13b}).

From academic point of view, the neutrosophic set and operators need
to be specified because is hard to be applied to the real
applications. So the concept of  interval neutrosophic sets
\cite{wan-05} which can represent uncertain, imprecise, incomplete
and inconsistent information was proposed. In this paper, we first
define interval neutrosophic soft sets (INS-sets) which is
generalizes the concept of the soft set, fuzzy soft set, interval
valued fuzzy soft set, intuitionistic fuzzy soft set, interval
valued intuitionistic fuzzy soft sets. Then, we introduce some
definitions and operations of interval neutrosophic soft sets. Some
properties of INS-sets which are connected to operations have been
established. Also, the aim of this paper is to investigate the
decision making based on interval valued neutrosophic soft sets. By
means of level soft sets, we develop an adjustable approach to
interval valued neutrosophic soft sets based decision making and a
examples are provided to illustrate the developed approach.

The relationship among interval neutrosophic soft set and other soft
sets is illustrated as;
$$
\begin{array}{rl}
\textrm{Soft set} & \subseteq\textrm {Fuzzy soft set}\\&  \subseteq
\textrm{Intuitionistic fuzzy soft set (Interval valued fuzzy soft
set)}
\\& \subseteq
\textrm{Interval valued intuitionistic fuzzy soft set }\\&
%\Rightarrow \textrm{Neutrosophic soft set}\\&
 \subseteq
\textrm{Interval valued neutrosophic soft set}
\end{array}
$$
Therefore, interval neutrosophic soft set is a generalization other
each the soft sets.
%-----------------------------------------------------------
\section{Preliminary}\label{ss}
%-----------------------------------------------------------
In this section, we present the basic definitions of neutrosophic
 set theory \cite{wan-05}, interval neutrosophic set theory \cite{wan-05} and
soft set theory \cite{mol-99} that are useful for subsequent
discussions. More detailed explanations related to this subsection
may be found in \cite{bro-13,bro-13b, cag-09c, jia-10,maj-13,maj-12,
wan-05}.
%-----------------------------------------------------------
\begin{defn} \cite{sma-05}
Let U be a space of points (objects), with a generic element in U
denoted by u. A neutrosophic sets(N-sets) A in U is characterized by
a truth-membership function $T_A$, a indeterminacy-membership
function $I_A$ and a falsity-membership function $F_A$. $T_A(u)$;
$I_A(x)$ and $F_A(u)$ are real standard or nonstandard subsets of
$[0,1]$.

There is no restriction on the sum of $T_A(u)$; $I_A(u)$ and
$F_A(u)$, so $0\leq sup T_A(u) + sup I_A(u) + supF_A(u)\leq 3$.
\end{defn}
%---------------------------------------------------------------
\begin{defn} \cite{wan-05}
Let U be a space of points (objects), with a generic element in U
denoted by u. An interval value neutrosophic set (IVN-sets) A in U
is characterized by truth-membership function $T_A$, a
indeterminacy-membership function $I_A$ and a falsity-membership
function $F_A$. For each point $u \in U$; $T_A$, $I_A$ and $F_A
\subseteq [0,1]$.

Thus, a IVN-sets over $U$ can be represented by the set of
$$
A= \{\langle T_A(u), I_A(u), F_A(u)\rangle/u:u\in U\}
$$
Here, $(T_A(u), I_A(u), F_A(u))$ is called interval value
neutrosophic number for all $u \in U$ and all interval value
neutrosophic numbers over $U$ will be denoted by $IVN(U)$.
\end{defn}
%---------------------------------------------------------------
\begin{exmp}\label{ex-softa}
Assume that the universe of discourse  $U=\{u_1, u_2\}$  where $u_1$
and characterises the quality,  $u_2$ indicates the prices of the
objects. It may be further assumed that the values of $u_1$ and
$u_2$ are subset of $[ 0, 1 ]$ and they are obtained from a expert
person. The expert construct an interval value neutrosophic set the
characteristics of the objects according to by truth-membership
function $T_A$, a indeterminacy-membership function $I_A$ and a
falsity-membership function $F_A$ as follows;

$A=\{\langle [0.1,1.0], [0.1,0.4], [0.4,0.7]\rangle/u_1,\langle
[0.6,0.9], [0.8,1.0], [0.4,0.6] \rangle/u_2\}$
\end{exmp}
%---------------------------------------------------------------
\begin{defn}\cite{wan-05} Let $A$ a interval neutrosophic sets.
Then, for all $u \in U$,
\begin{enumerate}
\item $A$ is empty, denoted
$A=\widetilde{{\emptyset}}$, is defined by
$$
\widetilde{{\emptyset}}= \{ <[ 0, 0],[ 1, 1],[ 1,1]>/u:u\in U\}
$$

\item $A$ is universal, denoted
$A=\widetilde{{E}}$, is defined by
$$
\widetilde{{E}}= \{<[ 1, 1],[ 0, 0],[ 0,0]
>/u:u\in U\}
$$

\item The complement of $A$ is denoted by
$\overline{A}$ and is defined by
$$
\begin{array}{rl}

\overline{A}=& \{<[infF_A(u),supF_A(u)],[ 1-sup I_A(u),1-inf
I_A(u)],\\&[infT_A(u),supT_A(u)]
>/u:u\in U\}\\
\end{array}
$$
\end{enumerate}
\end{defn}
%---------------------------------------------------------------------------------
\begin{defn} \cite{wan-05}
An interval neutrosophic set $A$ is contained in the other interval
neutrosophic set $B$, $A\widetilde{\subseteq }B$, if and only if
$$
\begin{array}{rl}
  \begin{array}{rl}
  &inf T_A(u) \leq inf T_B(u)\\&
sup T_A(u)\leq sup T_B(u)\\& inf
I_A(u)\geq inf I_B(u)\\
\end{array} &
\begin{array}{rl}
 &   sup I_A(u)\geq sup I_B(u)\\& inf F_A(u)\geq
inf F_B(u)\\& supF_A(u)\geq supF_B(u)
\end{array}
\end{array}
$$
for all $u \in U$.
\end{defn}
%---------------------------------------------------------------------------------
\begin{defn}
An interval neutrosophic number $X=(T_X,I_X,F_X)$ is larger than the
other interval neutrosophic number $Y=(T_Y,I_Y,F_Y)$, denoted
$X\widehat{\leq }Y$, if and only if
$$
\begin{array}{rl}
  \begin{array}{rl}
  &inf T_X \leq inf T_Y\\&
sup T_X\leq sup T_Y\\& inf
I_X\geq inf I_Y\\
\end{array} &
\begin{array}{rl}
 &   sup I_X\geq sup I_Y\\& inf F_X\geq
inf F_Y\\& supF_X\geq supF_Y
\end{array}
\end{array}
$$
\end{defn}
%---------------------------------------------------------------------------------
\begin{defn}\cite{wan-05} Let $A$ and $B$ be two interval neutrosophic sets.
Then, for all $u \in U$, $a \in R^+$,
\begin{enumerate}
\item Intersection of $A$ and $B$,
denoted by $A\widetilde{\cap} B$, is defined by
$$
\begin{array}{rl}
A\widetilde{\cap} B= &\{<[ min(inf T_A(u), i nf T_B(u)),min(sup
T_A(u), sup T_B(u))],\\&[ max(inf I_A(u), inf I_B(u)), max(sup
IA(x), sup I_B(u)) ],\\&[ max(inf F_A(u), inf F_B(u)),
max(supF_A(u),sup F_B(u))]
>/u:u\in U\}
\end{array}
$$

\item Union of $A$ and $B$,
denoted by $A\widetilde{\cup} B$, is defined by
$$
\begin{array}{rl}
A\widetilde{\cup} B= &\{ <[ max(inf T_A(u), inf T_B(u)),max(sup
T_A(u), sup T_B(u))],\\&[min(inf I_A(u), inf I_B(u)), min(sup
I_A(u), sup I_B(u))],\\&[min(inf F_A(u), inf F_B(u)), min(supF_A(u),
sup F_B(u))]
>/u:u\in U\}
\end{array}
$$

\item Difference of $A$ and $B$,
denoted by $A\widetilde{\setminus } B$, is defined by
$$
\begin{array}{rl}
A\widetilde{\setminus } B= &\{ <[  min(inf T_A(u), inf F_B(u)),
min(sup T_A(u), sup F_B( x))],\\&[max(inf I_A(u), 1- sup
I_B(u)),max(sup I_A(u), 1- inf I_B(u))],\\&[max(inf F_A(u), inf
T_B(u)), max(supF_A(u), sup T_B(u))]
>/u:u\in U\}
\end{array}
$$

\item Addition of $A$ and $B$,
denoted by $A\widetilde{+ } B$, is defined by
$$
\begin{array}{rl}
A\widetilde{+ } B= &\{<[  min(inf T_A(u) + inf T_B(u), 1), min(sup
T_A(u) + sup T_B(u), 1)],\\&[min(inf I_A(u) + inf I_B(u), 1),
min(sup I_A(u) + sup I_B(u), 1)],\\&[min(inf F_A(u) + inf F_B(u),
1), min(supF_A(u) + supF_B(u), 1)]
>/u:u\in U\}
\end{array}
$$

\item Scalar multiplication of $A$,
denoted by $A\widetilde{. } a$, is defined by
$$
\begin{array}{rl}
A\widetilde{. }a=&\{ <[  min(inf T_A(u).a, 1), min(sup T_A(u). a,
1)],\\&[min(infI_A(u).a, 1), min(sup I_A(u).a, 1)],\\&[min(inf
F_A(u).a, 1),min(sup F_A(u) . a, 1)]
>/u:u\in U\}
\end{array}
$$

\item Scalar division of $A$,
denoted by $A\widetilde{/ }a$, is defined by
$$
\begin{array}{rl}
A\widetilde{/} a= &\{ <[  min(inf T_A(u)/a, 1), min(sup T_A(u)/ a,
1)],\\&[min(infI_A(u)/a, 1), min(sup I_A(u)/a, 1)],\\&[min(inf
F_A(u)/a, 1),min(sup F_A(u) / a, 1)]
>/u:u\in U\}
\end{array}
$$

\item Truth-Favorite of $A$,
denoted by $\widetilde{\triangle} A$, is defined by
$$
\begin{array}{rl}
\widetilde{\triangle} A= &\{ <[  min(inf T_A(u) + inf I_A(u), 1),
min(sup T_A(u) + sup I_A(u), 1)],[0, 0],\\&[inf F_A(u) , supF_A(u)]
>/u:u\in U\}
\end{array}
$$

\item False-Favorite of $A$,
denoted by $\widetilde{\nabla} A$, is defined by
$$
\begin{array}{rl}
\widetilde{\nabla} A = &\{ <[   inf T_A(u) ,sup T_A(u) ],[0,
0],\\&[min(inf F_A(u) + inf I_A(u), 1),min(supF_A(u) + sup I_A(u),
1)]
>/u:u\in U\}
\end{array}
$$
\end{enumerate}
\end{defn}
%---------------------------------------------------------------
\begin{defn}\cite{mol-99}
Let $U$ be an initial universe,  $P(U)$ be the power set of $U$, $E$
be a set of all parameters and $X\subseteq E$. Then a soft set $F_X$
over $U$ is a set defined by a function representing a mapping
$$f_X: E\to P(U)\, such \,that \,f_X(x)=\emptyset \,\,if\,\, x\notin
X$$ Here, $f_X$ is called approximate function of the soft set
$F_X$, and the value $f_X(x)$ is a set called \emph{$x$-element} of
the soft set for all  $x \in E$. It is worth noting that the sets is
worth noting that the sets $f_X(x)$ may be arbitrary. Some of them
may be empty, some may have nonempty intersection. Thus, a soft set
over $U$ can be represented by the set of ordered pairs

$$ F_X= \{(x, f_X(x)): x\in E, f_X(x)\in P(U)\}$$
\end{defn}
%---------------------------------------------------------------
\begin{exmp}\label{ex-soft1}
Suppose that $U=\{u_1,u_2,u_3,u_4,u_5,u_6\}$ is the universe
contains six house under consideration in an real agent and
$E=\{x_1=cheap, x_2=beatiful, x_3=green surroundings, x_4=costly,
x_5= large\}$.

A customer to select a house from the real agent. then, he can
construct a soft set $F_X$ that describes the characteristic of
houses according to own requests. Assume that $f_X(x_1)=
\{u_1,u_2\}$, $f_X(x_2)= \{u_1 \}$, $ f_X(x_3)= \emptyset$, $
f_X(x_4)= U$, $ \{u_1,u_2,u_3,u_4,u_5\}$ then the soft-set $F_X$ is
written by
$$
F_X=\{(x_1, \{u_1,u_2\}),(x_2, \{u_1,u_4,u_5,u_6\}), (x_4,U),
(x_5,\{u_1,u_2,u_3,u_4,u_5\})\}
$$

The tabular representation of the soft set $F_X$ is as follow:
$$
\begin{tabular}{|c|cccccc|}
  \hline
 \(U\) &  \(u_1\) &
\(u_2\)&  \(u_3\) & \(u_4\)&  \(u_5\) & \(u_6\)\\
    \hline \(x_1\)  &
\(1\)& \(1\)& \(0\)& \(0\)& \(0\)& \(0\) \\
    \(x_2\) &
\(1\)& \(0\)& \(0\)& \(1\)& \(1\)& \(1\) \\
    \(x_3\)  &
\(0\)& \(0\)& \(0\)& \(0\)& \(0\)& \(0\)\\
\(x_4\)  &
\(1\)& \(1\)& \(1\)& \(1\)& \(1\)& \(1\)\\
  \(x_5\)  &
\(1\)& \(1\)& \(1\)& \(1\)& \(1\)& \(0\)\\
   \hline
    \end{tabular}
   $$
\begin{center}
\footnotesize{\emph{\emph{Table 1}: The tabular representation of
the soft set $F_X$}}
\end{center}
\end{exmp}

\begin{defn}\cite{maj-02}\label{cho}
Let $U=\{u_1,u_2,...,u_k\}$ be an initial universe of objects,
$E=\{x_1, x_2,... ,x_m\}$ be a set of parameters and $F_X$ be a soft
set over U. For any $x_j \in E$, $f_X(x_j)$ is a subset of U. Then,
the choice value of an object $u_i\in U$ is $c_i$, given by $c_i =
\sum_j u_{ij}$, where $u_{ij}$ are the entries in the table of the
reduct-soft-set. That is,
$$
u_{ij}=\left\{
\begin{array}{ll}
1, & u_i\in f_X(x_j)\\
0, & u_i \notin f_X(x_j)
\end{array}\right.
$$
\end{defn}
%---------------------------------------------------------------
\begin{exmp}\label{ex-soft2}
Consider the above Example \ref{ex-soft1}. Clearly,

$$c_1= \sum_{j=1}^5 u_{1j}=4,$$ $$ c_3=c_6= \sum_{j=1}^5 u_{3j}=
\sum_{j=1}^5 u_{6j}=2,$$
$$
   c_2=c_4=c_5= \sum_{j=1}^5 u_{2j}=\sum_{j=1}^5 u_{4j}= \sum_{j=1}^5 u_{5j}=3$$.
\end{exmp}
%----------------------------------------------------------
\begin{defn}\cite{cag-09a} Let $F_X$ and $F_Y$ be two soft sets. Then,
\begin{enumerate} \item
Complement of $F_X$ is denoted by $F_X^{\tilde{c}}$. Its approximate
function $f_{X^c}(x)= U\setminus f_X(x)\quad\textrm{for all } x \in
E $
\item
Union of $F_X$ and $F_Y$ is denoted by $F_X \tilde{\cup} F_Y$. Its
approximate function $f_{X\tilde{\cup} F_Y}$ is defined by
$$
f_{X \tilde{\cup} Y}(x)= f_X(x)\cup f_Y(x)\quad\textrm{for all } x
\in E.
$$

\item Intersection of $F_X$ and $F_Y$ is
denoted by $F_X\tilde{\cap}F_Y$. Its approximate function $f_{X
\widetilde{\cap} Y}$ is defined by
$$
f_{X \tilde{\cap} Y}(x)= f_X(x)\cap f_Y(x)\quad\textrm{for all } x
\in E.
$$
\end{enumerate}
\end{defn}
%---------------------------------------------------------------
\section{Interval-valued neutrosophic soft sets}
In this section, we give interval valued neutrosophic soft sets
($ivn-$soft sets) which is a combination of an interval valued
neutrosophic sets\cite{wan-05}  and a soft sets\cite{mol-99}. Then,
we introduce some definitions and operations of $ivn-$soft sets
sets. Some properties of $ivn-$soft sets which are connected to
operations have been established. Some of it is quoted from
\cite{bro-13,bro-13b, cag-09a, cag-09c, jia-10,maj-13,maj-12,
wan-05}.
%---------------------------------------------------------------
\begin{defn}
Let $U$ be an initial universe set, $IVN( U )$ denotes the set of
all interval valued neutrosophic sets of $U$ and $E$ be a set of
parameters that are describe the elements of $U$. An interval valued
neutrosophic soft sets($ivn$-soft  sets) over $U$ is a set defined
by a set valued function $\Upsilon_K$ representing a mapping
$$
\upsilon_K: E\to IVN(U)
$$
It can be written a set of ordered pairs
$$
\Upsilon_K= \{(x, \upsilon_K(x)): x\in E\}
$$
Here, $\upsilon_K$, which is interval valued neutrosophic sets, is
called approximate function of the $ivn$-soft  sets $\Upsilon_K$ and
$\upsilon_K(x)$ is called $x$-approximate value of $x \in E$. The
subscript $K$ in the $\upsilon_K$ indicates that $\upsilon_K$ is the
approximate function of $\Upsilon_K$.

Generally, $\upsilon_K$, $\upsilon_L$, $\upsilon_M$,... will be used
as an approximate functions of $\Upsilon_K$, $\Upsilon_L$,
$\Upsilon_M$, ..., respectively.
\end{defn}
Note that the sets of all $ivn$-soft  sets over $U$ will be denoted
by $IVNS(U) $.
%---------------------------------------------------------------

Now let us give the following example for  $ivn$-soft  sets.
%---------------------------------------------------------------
\begin{exmp}\label{ex-soft}
Let $U=\{u_1,u_2\}$ be set of houses under consideration and $E$ is
a set of parameters which is a neutrosophic word. Consider
$E=\{x_1=cheap, x_2=beatiful, x_3=green surroundings, x_4=costly,
x_5= large\}$. In this case,  we give an ($ivn$-soft   sets)
$\Upsilon_K$ over $U$ as;
$$
\begin{array}{lr}
\Upsilon_K = &\{(x_1,\{\langle [0.6,0.8], [0.8,0.9], [0.1,0.5]
\rangle /u_1,\langle [0.5,0.8], [0.2,0.9], [0.1,0.7] \rangle /u_2
\}),\\&(x_2,\{\langle [0.1,0.4], [0.5,0.8], [0.3,0.7] \rangle
/u_1,\langle [0.1,0.9], [0.6,0.9], [0.2,0.3] \rangle /u_2
\}),\\&(x_3,\{\langle [0.2,0.9], [0.1,0.5], [0.7,0.8]  \rangle
/u_1,\langle [0.4,0.9], [0.1,0.6], [0.5,0.7] \rangle
/u_2\}),\\&(x_4,\{\langle [0.6,0.9], [0.6,0.9], [0.6,0.9] \rangle
/u_1,\langle[0.5,0.9], [0.6,0.8], [0.1,0.8] \rangle
/u_2\}),\\&(x_5,\{\langle [0.0,0.9], [1.0,1.0], [1.0,1.1] \rangle
/u_1,\langle [0.0,0.9], [0.8,1.0], [0.2,0.5] \rangle /u_2\})\}
\end{array}
$$
The tabular representation of the $ivn$-soft set $\Upsilon_K$ is as
follow:
$$
\begin{tabular}{|c|cc|}
  \hline
 \(U\) &  \(u_1\) &
\(u_2\)\\
    \hline \(x_1\)  &
\(\langle [0.6,0.8], [0.8,0.9], [0.1,0.5]\rangle\)& \(\langle
[0.5,0.8], [0.2,0.9], [0.1,0.7]
\rangle\) \\
    \(x_2\) &
\(\langle [0.1,0.4], [0.5,0.8], [0.3,0.7] \rangle\)& \(\langle
[0.1,0.9], [0.6,0.9], [0.2,0.3]
\rangle\) \\
    \(x_3\)  &
\(\langle [0.2,0.9], [0.1,0.5], [0.7,0.8] \rangle\)& \(\langle
[0.4,0.9], [0.1,0.6], [0.5,0.7]
\rangle\) \\
    \(x_4\)  &
\(\langle [0.6,0.9], [0.6,0.9], [0.6,0.9] \rangle\)& \(\langle
[0.5,0.9], [0.6,0.8], [0.1,0.8]
\rangle\) \\
   \(x_5\) &
\(\langle [0.0,0.9], [1.0,1.0], [1.0,1.1]\rangle\)& \(\langle
 [0.0,0.9], [0.8,1.0], [0.2,0.5]
\rangle\) \\
     \hline
\end{tabular}
$$
\begin{center}
\footnotesize{\emph{\emph{Table 1}: The tabular representation of
the $ivn$-soft set $\Upsilon_K$}}
\end{center}
\end{exmp}
%---------------------------------------------------------------
\begin{defn}
Let $\Upsilon_K \in IVNS(U)$. If
$\upsilon_K(x)=\widetilde{{\emptyset}}$ for all $x\in E$, then $N$
is called an empty $ivn$-soft  set, denoted by
${\widehat{\emptyset}}$.
\end{defn}
%-----------------------------------------------------------
\begin{defn}
Let $\Upsilon_K \in IVNS(U)$. If  $\upsilon_K(x)=\widetilde{{E}}$
for all $x\in E$, then $\Upsilon_K$ is called a universal $ivn$-soft
set, denoted by $\Upsilon_{\Upsilon_{\widehat{E}}}$.
\end{defn}
%---------------------example-------------------------------
\begin{exmp}\label{ex-soft}
Assume that $U=\{u_1, u_2\}$ is a universal set and $E=\{x_1, x_2,
x_3,$ $x_4,x_5\}$ is a set of all parameters. Consider the tabular
representation of the $\Upsilon_{\widehat{\emptyset}}$ is as
follows;
$$
\begin{tabular}{|c|cc|}
  \hline
 \(U\) &  \(u_1\) &
\(u_2\)\\
    \hline \(x_1\)  &
\(\langle [0.0,0.0], [1.0,1.0], [1.0,1.0]\rangle\)& \(\langle
[0.0,0.0], [1.0,1.0], [1.0,1.0]
\rangle\) \\
    \(x_2\) &
\(\langle [0.0,0.0], [1.0,1.0], [1.0,1.0] \rangle\)& \(\langle
[0.0,0.0], [1.0,1.0], [1.0,1.0]
\rangle\) \\
    \(x_3\)  &
\(\langle [0.0,0.0], [1.0,1.0], [1.0,1.0]\rangle\)& \(\langle
[0.0,0.0], [1.0,1.0], [1.0,1.0]
\rangle\) \\
    \(x_4\)  &
\(\langle[0.0,0.0], [1.0,1.0], [1.0,1.0] \rangle\)& \(\langle
[0.0,0.0], [1.0,1.0], [1.0,1.0]
\rangle\) \\
   \(x_5\) &
\(\langle [0.0,0.0], [1.0,1.0], [1.0,1.0]\rangle\)& \(\langle
[0.0,0.0], [1.0,1.0], [1.0,1.0]
\rangle\) \\
    \hline
\end{tabular}
$$
\begin{center}
\footnotesize{\emph{\emph{Table 2}: The tabular representation of
the $ivn$-soft set $\Upsilon_{\widehat{\emptyset}}$}}
\end{center}
The tabular representation of the $\Upsilon_{\widehat{E}}$ is as
follows;
$$
\begin{tabular}{|c|cc|}
  \hline
 \(U\) &  \(u_1\) &
\(u_2\)\\
    \hline \(x_1\)  &
\(\langle [1.0,1.0], [0.0,0.0], [0.0,0.0] \rangle\)& \(\langle
[1.0,1.0], [0.0,0.0], [0.0,0.0]
\rangle\) \\
    \(x_2\) &
\(\langle [1.0,1.0], [0.0,0.0], [0.0,0.0] \rangle\)& \(\langle
[1.0,1.0], [0.0,0.0], [0.0,0.0]
\rangle\) \\
    \(x_3\)  &
\(\langle [1.0,1.0], [0.0,0.0], [0.0,0.0] \rangle\)& \(\langle
[1.0,1.0], [0.0,0.0], [0.0,0.0]
\rangle\) \\
    \(x_4\)  &
\(\langle [1.0,1.0], [0.0,0.0], [0.0,0.0]\rangle\)& \(\langle
[1.0,1.0], [0.0,0.0], [0.0,0.0]
\rangle\) \\
   \(x_5\) &
\(\langle [1.0,1.0], [0.0,0.0],[0.0,0.0]\rangle\)& \(\langle
[1.0,1.0], [0.0,0.0], [0.0,0.0]
\rangle\) \\
    \hline
\end{tabular}
$$
\begin{center}
\footnotesize{\emph{\emph{Table 3}: The tabular representation of
the $ivn$-soft set $\Upsilon_{\widehat{E}}$}}
\end{center}
\end{exmp}
%-----------------------------------------------------------
\begin{defn}\label{subset}
Let $\Upsilon_K, \Upsilon_L\in IVNS(U)$. Then, $\Upsilon_K$ is an
$ivn$-soft subset of $\Upsilon_L $, denoted by $\Upsilon_K
\widehat{{\subseteq}} \Upsilon_L$, if
$\upsilon_K(x)\widetilde{\subseteq} \upsilon_L(x)$ for all $x \in
E$.
\end{defn}
%---------------------------------------------------------------
\begin{exmp}\label{ex-soft}
Assume that $U=\{u_1, u_2\}$ is a universal set and $E=\{x_1, x_2,
x_3,$ $x_4,x_5\}$ is a set of all parameters. Consider the tabular
representation of the $\Upsilon_K$ is as follows;
$$
\begin{tabular}{|c|cc|}
  \hline
 \(U\) &  \(u_1\) &
\(u_2\)\\
    \hline \(x_1\)  &
\(\langle [0.5,0.7], [0.8,0.9], [0.2,0.5]\rangle\)& \(\langle
[0.3,0.6], [0.3,0.9], [0.2,0.8]
\rangle\) \\
    \(x_2\) &
\(\langle [0.0,0.3], [0.6,0.8], [0.3,0.9] \rangle\)& \(\langle
[0.1,0.8], [0.8,0.9], [0.3,0.5]
\rangle\) \\
    \(x_3\)  &
\(\langle [0.1,0.7], [0.4,0.5], [0.8,0.9] \rangle\)& \(\langle
[0.2,0.5], [0.5,0.7], [0.6,0.8]
\rangle\) \\
    \(x_4\)  &
\(\langle [0.2,0.4], [0.7,0.9], [0.6,0.9] \rangle\)& \(\langle
[0.3,0.9], [0.6,0.9], [0.3,0.9]
\rangle\) \\
   \(x_5\) &
\(\langle [0.0,0.2], [1.0,1.0], [1.0,1.0]\rangle\)& \(\langle
 [0.0,0.1], [0.9,1.0], [0.2,0.9]
\rangle\) \\
     \hline
\end{tabular}
$$
\begin{center}
\footnotesize{\emph{\emph{Table 4}: The tabular representation of
the $ivn$-soft set $\Upsilon_{K}$}}
\end{center}
The tabular representation of the $\Upsilon_L$ is as follows;
$$
\begin{tabular}{|c|cc|}
  \hline
 \(U\) &  \(u_1\) &
\(u_2\)\\
    \hline \(x_1\)  &
\(\langle [0.6,0.8], [0.8,0.9], [0.1,0.5]\rangle\)& \(\langle
[0.5,0.8], [0.2,0.9], [0.1,0.7]
\rangle\) \\
    \(x_2\) &
\(\langle [0.1,0.4], [0.5,0.8], [0.3,0.7] \rangle\)& \(\langle
[0.1,0.9], [0.6,0.9], [0.2,0.3]
\rangle\) \\
    \(x_3\)  &
\(\langle [0.2,0.9], [0.1,0.5], [0.7,0.8] \rangle\)& \(\langle
[0.4,0.9], [0.1,0.6], [0.5,0.7]
\rangle\) \\
    \(x_4\)  &
\(\langle [0.6,0.9], [0.6,0.9], [0.6,0.9] \rangle\)& \(\langle
[0.5,0.9], [0.6,0.8], [0.1,0.8]
\rangle\) \\
   \(x_5\) &
\(\langle [0.0,0.9], [1.0,1.0], [1.0,1.0]\rangle\)& \(\langle
 [0.0,0.9], [0.8,1.0], [0.2,0.5]
\rangle\) \\
     \hline
\end{tabular}
$$
\begin{center}
\footnotesize{\emph{\emph{Table 5}: The tabular representation of
the $ivn$-soft set $\Upsilon_{L}$}}
\end{center}
Clearly, by Definition \ref{subset}, we have $\Upsilon_K
\widehat{{\subseteq}} \Upsilon_L$.
\end{exmp}
%-----------------------------------------------------------
\begin{rem}
$\Upsilon_K \widehat{{\subseteq}} \Upsilon_L$ does not imply that
every element of $\Upsilon_K$ is an element of $\Upsilon_L$ as in
the definition of the classical subset.
\end{rem}
%-----------------------------------------------------------
\begin{prop}Let $\Upsilon_K, \Upsilon_L,\Upsilon_M \in IVNS(U)$. Then,
\begin{enumerate}
\item  $\Upsilon_K \widehat{\subseteq} {\Upsilon_{\widehat{E}}}$
\item  $\Upsilon_{\widehat{\emptyset}} \widehat{\subseteq} \Upsilon_K $
\item  $\Upsilon_K\widehat {\subseteq} \Upsilon_K$
\item  $\Upsilon_K \widehat{\subseteq} \Upsilon_L$ and $\Upsilon_L \widehat{\subseteq} \Upsilon_M
\Rightarrow \Upsilon_K \widehat{\subseteq} \Upsilon_M$
\end{enumerate}
\end{prop}
%---------------------------------------------------------------
\begin{pf}They can be proved easily by using the approximate function of the $ivn$-soft  sets.
\end{pf}
%-----------------------------------------------------------
\begin{defn}
Let $\Upsilon_K, \Upsilon_L \in IVNS(U)$. Then, $\Upsilon_K$ and
$\Upsilon_L$ are $ivn$-soft  equal, written as $\Upsilon_K
{=}\Upsilon_L$, if and only if $\upsilon_K(x)= \upsilon_L(x)$ for
all $x \in E$.
\end{defn}
%-----------------------------------------------------------
\begin{prop}Let $\Upsilon_K, \Upsilon_L,\Upsilon_M \in IVNS(U)$. Then,
\begin{enumerate}
\item $\Upsilon_K {=}\Upsilon_L$ and $\Upsilon_L =\Upsilon_M  \Leftrightarrow \Upsilon_K=\Upsilon_M$
\item $\Upsilon_K \widehat{\subseteq} \Upsilon_L$ and $\Upsilon_L \widehat{\subseteq} \Upsilon_K  \Leftrightarrow \Upsilon_K =\Upsilon_L$
\end{enumerate}
\end{prop}
\begin{pf}The proofs are trivial.
\end{pf}
%-----------------------------------------------------------
\begin{defn}\label{comp}
Let $\Upsilon_K \in IVNS(U)$. Then, the complement
$\Upsilon_K^{\widehat{{c}}}$ of $\Upsilon_K$ is an $ivn$-soft  set
such that
$$
\upsilon_K^{\widehat{{c}}}(x)=\overline{\upsilon_K}(x),\textrm{ for
all } x \in E.
$$
\end{defn}
%---------------------------------------------------------------
\begin{exmp}\label{ex-soft-2}
Consider the above Example \ref{ex-soft}, the complement
$\Upsilon_L^{\widehat{{c}}}$ of $\Upsilon_L$ can be represented into
the following table;
$$
\begin{tabular}{|c|cc|}
  \hline
 \(U\) &  \(u_1\) &
\(u_2\)\\
    \hline \(x_1\)  &
\(\langle [0.1,0.5], [0.1,0.2], [0.6,0.8]\rangle\)& \(\langle
[0.1,0.7], [0.1,0.8], [0.5,0.8]
\rangle\) \\
    \(x_2\) &
\(\langle [0.3,0.7], [0.2,0.5], [0.1,0.4] \rangle\)& \(\langle
[0.2,0.3], [0.1,0.4], [0.1,0.9]
\rangle\) \\
    \(x_3\)  &
\(\langle[0.7,0.8] , [0.5,0.9], [0.2,0.9] \rangle\)& \(\langle
[0.5,0.7], [0.4,0.9],[0.4,0.9]
\rangle\) \\
    \(x_4\)  &
\(\langle [0.6,0.9], [0.1,0.4], [0.6,0.9] \rangle\)& \(\langle
[0.1,0.8], [0.2,0.4],[0.5,0.9]
\rangle\) \\
   \(x_5\) &
\(\langle [1.0,1.0], [0.0,0.0],[0.0,0.9] \rangle\)& \(\langle
 [0.2,0.5], [0.0,0.2],[0.0,0.9]
\rangle\) \\
     \hline
\end{tabular}
$$
\begin{center}
\footnotesize{\emph{\emph{Table 6}: The tabular representation of
the $ivn$-soft set $\Upsilon_L^{\widehat{{c}}}$}}
\end{center}
\end{exmp}
%-----------------------------------------------------------
\begin{prop}Let $\Upsilon_K \in IVNS(U)$. Then,
\begin{enumerate}
\item  $(\Upsilon_K^{\widehat{{c}}})^{\widehat{{c}}}= \Upsilon_K$
\item  $ \Upsilon_{\widehat{\emptyset}} ^{\widehat{{c}}} = \Upsilon_{\widehat{E}} $
\item  $ \Upsilon_{\widehat{E}} ^{\widehat{{c}}} = \Upsilon_{\widehat{\emptyset}} $
\end{enumerate}
\end{prop}\begin{pf}By using the fuzzy approximate functions of the $ivn$-soft  set,
the proofs can be straightforward.
\end{pf}
%-----------------------------------------------------------
\begin{thm}Let $\Upsilon_K \in IVNS(U)$. Then,
 $\Upsilon_K\widehat{\subseteq }\Upsilon_L \Leftrightarrow
 \Upsilon_L^{\widehat{{c}}}\widehat{\subseteq }\Upsilon_K^{\widehat{{c}}} $
\end{thm}
\begin{pf}By using the fuzzy approximate functions of the $ivn$-soft  set,
the proofs can be straightforward.
\end{pf}
%-----------------------------------------------------------
\begin{defn}\label{union}
Let $\Upsilon_K, \Upsilon_L \in IVNS(U)$. Then, union of
$\Upsilon_K$ and $\Upsilon_L$, denoted $ \Upsilon_K \widehat{\cup}
\Upsilon_L$, is defined by
$$
\upsilon_{K \widehat{\cup}
L}(x)=\upsilon_{K}(x)\widetilde{\cup}\upsilon_{L}(x)\quad\textrm{for
all } x \in E.
$$
\end{defn}
%---------------------------------------------------------------
\begin{exmp}\label{ex-soft-3}
Consider the above Example \ref{ex-soft}, the union of $\Upsilon_K$
and $\Upsilon_L$, denoted $ \Upsilon_K \widehat{\cup} \Upsilon_L$,
can be represented into the following table;
$$
\begin{tabular}{|c|cc|}
  \hline
 \(U\) &  \(u_1\) &
\(u_2\)\\
    \hline \(x_1\)  &
\(\langle [0.6,0.8], [0.8,0.9], [0.1,0.5]\rangle\)& \(\langle
[0.5,0.8], [0.2,0.9], [0.1,0.7]
\rangle\) \\
    \(x_2\) &
\(\langle [0.1,0.4], [0.5,0.8], [0.3,0.7] \rangle\)& \(\langle
[0.1,0.9], [0.6,0.9], [0.2,0.3]
\rangle\) \\
    \(x_3\)  &
\(\langle [0.2,0.9], [0.1,0.5], [0.7,0.8] \rangle\)& \(\langle
[0.4,0.9], [0.1,0.6], [0.5,0.7]
\rangle\) \\
    \(x_4\)  &
\(\langle [0.6,0.9], [0.6,0.9], [0.6,0.9] \rangle\)& \(\langle
[0.5,0.9], [0.6,0.8], [0.1,0.8]
\rangle\) \\
   \(x_5\) &
\(\langle [0.0,0.9], [1.0,1.0], [1.0,1.0]\rangle\)& \(\langle
 [0.0,0.9], [0.8,1.0], [0.2,0.5]
\rangle\) \\
     \hline
\end{tabular}
$$
\begin{center}
\footnotesize{\emph{\emph{Table 7}: The tabular representation of
the $ivn$-soft set $\Upsilon_K \widehat{\cup} \Upsilon_L$}}
\end{center}
\end{exmp}
%-----------------------------------------------------------
\begin{thm}Let $\Upsilon_K, \Upsilon_L \in IVNS(U)$. Then,
$\Upsilon_K  \widehat{\cup} \Upsilon_L $ is the smallest $ivn-$ soft
set containing both $\Upsilon_K$ and $\Upsilon_L$.
\end{thm}
\begin{pf}The proofs can be easily obtained from Definition \ref{union}.
\end{pf}
%-----------------------------------------------------------
\begin{prop}Let $\Upsilon_K, \Upsilon_L,\Upsilon_M \in IVNS(U)$. Then,
\begin{enumerate}
\item  $ \Upsilon_K\widehat{\cup} \Upsilon_K = \Upsilon_K $
\item  $ \Upsilon_K \widehat{\cup} \Upsilon_{\widehat{\emptyset}} = \Upsilon_K$
\item  $ \Upsilon_K \widehat{\cup} {\Upsilon_{\widehat{E}}} = {\Upsilon_{\widehat{E}}} $
\item  $ \Upsilon_K \widehat{\cup} \Upsilon_L= \Upsilon_L\widehat{\cup} \Upsilon_K $
\item  $(\Upsilon_K \widehat{\cup} \Upsilon_L)\widehat{\cup} \Upsilon_M= \Upsilon_K \widehat{\cup} (\Upsilon_L\widehat{\cup} \Upsilon_M) $
\end{enumerate}
\end{prop}\begin{pf}The proofs can be easily obtained from Definition \ref{union}.
\end{pf}%-----------------------------------------------------------
\begin{defn}\label{intersection}
Let $\Upsilon_K, \Upsilon_L \in IVNS(U)$. Then, intersection of
$\Upsilon_K$ and $\Upsilon_L$, denoted $ \Upsilon_K \widehat{\cap}
\Upsilon_L$, is defined by
$$
\upsilon_{K \widehat{\cap}
L}(x)=\upsilon_K(x)\widetilde{\cap}\upsilon_{L}(x)\quad\textrm{for
all } x \in E.
$$
\end{defn}
%---------------------------------------------------------------
\begin{exmp}\label{ex-soft-3}
Consider the above Example \ref{ex-soft}, the intersection of
$\Upsilon_K$ and $\Upsilon_L$, denoted $ \Upsilon_K \widehat{\cap}
\Upsilon_L$, can be represented into the following table;
$$
\begin{tabular}{|c|cc|}
  \hline
 \(U\) &  \(u_1\) &
\(u_2\)\\
    \hline \(x_1\)  &
\(\langle [0.5,0.7], [0.8,0.9], [0.2,0.5]\rangle\)& \(\langle
[0.3,0.6], [0.3,0.9], [0.2,0.8]
\rangle\) \\
    \(x_2\) &
\(\langle [0.0,0.3], [0.6,0.8], [0.3,0.9] \rangle\)& \(\langle
[0.1,0.8], [0.8,0.9], [0.3,0.5]
\rangle\) \\
    \(x_3\)  &
\(\langle [0.1,0.7], [0.4,0.5], [0.8,0.9] \rangle\)& \(\langle
[0.2,0.5], [0.5,0.7], [0.6,0.8]
\rangle\) \\
    \(x_4\)  &
\(\langle [0.2,0.4], [0.7,0.9], [0.6,0.9] \rangle\)& \(\langle
[0.3,0.9], [0.6,0.9], [0.3,0.9]
\rangle\) \\
   \(x_5\) &
\(\langle [0.0,0.2], [1.0,1.0], [1.0,1.0]\rangle\)& \(\langle
 [0.0,0.1], [0.9,1.0], [0.2,0.9]
\rangle\) \\
     \hline
\end{tabular}
$$
\begin{center}
\footnotesize{\emph{\emph{Table 8}: The tabular representation of
the $ivn$-soft set $\Upsilon_K \widehat{\cap} \Upsilon_L$}}
\end{center}
\end{exmp}
%-----------------------------------------------------------
\begin{prop}Let $\Upsilon_K, \Upsilon_L \in IVNS(U)$. Then,
$\Upsilon_K  \widehat{\cap} \Upsilon_L $ is the largest $ivn-$ soft
set containing both $\Upsilon_K$ and $\Upsilon_L$.
\end{prop}
\begin{pf}The proofs can be easily obtained from Definition \ref{intersection}.
\end{pf}
%-----------------------------------------------------------
\begin{prop}Let $\Upsilon_K, \Upsilon_L,\Upsilon_M \in IVNS(U)$. Then,
\begin{enumerate}
\item  $ \Upsilon_K \widehat{\cap} \Upsilon_K = \Upsilon_K $
\item  $ \Upsilon_K \widehat{\cap} \Upsilon_{\widehat{\emptyset}} = \Upsilon_{\widehat{\emptyset}}$
\item  $ \Upsilon_K \widehat{\cap} {\Upsilon_{\widehat{E}}} = \Upsilon_K  $
\item  $ \Upsilon_K \widehat{\cap} \Upsilon_L= \Upsilon_L\widehat{\cap}\Upsilon_K $
\item  $(\Upsilon_K \widehat{\cap} \Upsilon_L)\widehat{\cap} \Upsilon_M=\Upsilon_K \widehat{\cap} (\Upsilon_L\widehat{\cap} \Upsilon_M) $
\end{enumerate}
\end{prop}
\begin{pf} The proof of the Propositions 1- 5 are obvious.
\end{pf}
%-----------------------------------------------------------
\begin{rem}
Let $\Upsilon_K \in IVNS(U)$. If $\Upsilon_K \neq
\Upsilon_{\widehat{\emptyset}} $ or $\Upsilon_K \neq
\Upsilon_{\widehat {E}}$, then $\Upsilon_K \widehat{\cup}
\Upsilon_K^{\widehat{c}}\neq \Upsilon_{\widehat{E}}$ and $\Upsilon_K
\widehat{\cap}\Upsilon_K^{\widehat{c}}\neq
\Upsilon_{\widehat{\emptyset}}$.
\end{rem}
%-----------------------------------------------------------
\begin{prop}Let  $\Upsilon_K, \Upsilon_L \in IVNS(U)$. Then, De Morgan's laws are valid
\begin{enumerate}
\item $ (\Upsilon_K \widehat{\cup} \Upsilon_L)^{\widehat{c}}=
\Upsilon_K^{\widehat{c}} \widehat{\cap} \Upsilon_L^{\widehat{c}} $
\item $ (\Upsilon_K \widehat{\cap} \Upsilon_L)^{\widehat{c}} =
\Upsilon_K^{\widehat{c}} \widehat{\cup} \Upsilon_L^{\widehat{c}}.$
\end{enumerate}
\end{prop}
\begin{pf}The proofs can be easily obtained from Definition
\ref{comp}, Definition \ref{union} and Definition
\ref{intersection}.
\end{pf}
%-----------------------------------------------------------
\begin{prop}Let $\Upsilon_K, \Upsilon_L,\Upsilon_M \in IVNS(U)$. Then,
\begin{enumerate}
\item$ \Upsilon_K \widehat{\cup} (\Upsilon_L \widehat{\cap} \Upsilon_M)=
(\Upsilon_K \widehat{\cup} \Upsilon_L) \widehat{\cap} (\Upsilon_K
\widehat{\cup} \Upsilon_M) $
\item$ \Upsilon_K \widehat{\cap} (\Upsilon_L \widehat{\cup}\Upsilon_M)=
(\Upsilon_K\widehat{\cap} \Upsilon_L)\widehat{\cup} (\Upsilon_K
\widehat{\cap} \Upsilon_M) $
\item$ \Upsilon_K \widehat{\cup} (\Upsilon_K \widehat{\cap} \Upsilon_L)=
\Upsilon_K$
\item$ \Upsilon_K \widehat{\cap} (\Upsilon_K
\widehat{\cup}\Upsilon_L)= \Upsilon_K$
\end{enumerate}
\end{prop}
\begin{pf}The proofs can be easily obtained from Definition \ref{union} and Definition \ref{intersection}.
\end{pf}
%---------------------------------------------------------
\begin{defn}\label{or}
Let $\Upsilon_K, \Upsilon_L \in IVNS(U)$. Then, OR operator of
$\Upsilon_K$ and $\Upsilon_L$, denoted $ \Upsilon_K
\widehat{\bigvee} \Upsilon_L$, is defined by a set valued function
$\Upsilon_O$ representing a mapping
$$
\upsilon_O: E\times E\to IVN(U)
$$
where
$$
\upsilon_O(x,y)=\upsilon_K(x)\widetilde{\cup}\upsilon_{L}(y)\quad\textrm{for
all }(x,y) \in E\times E.
$$
\end{defn}
%---------------------------------------------------------
\begin{defn}\label{and}
Let $\Upsilon_K, \Upsilon_L \in IVNS(U)$. Then, AND operator of
$\Upsilon_K$ and $\Upsilon_L$, denoted $ \Upsilon_K
\widehat{\bigwedge} \Upsilon_L$, is defined by is defined by a set
valued function $\Upsilon_A$ representing a mapping
$$
\upsilon_A: E\times E\to IVN(U)
$$
where
$$
\upsilon_A(x,y)=\upsilon_K(x)\widetilde{\cap}\upsilon_{L}(y)\quad\textrm{for
all }(x,y) \in E\times E.
$$
\end{defn}
%-----------------------------------------------------------
\begin{prop}Let $\Upsilon_K, \Upsilon_L,\Upsilon_M \in IVNS(U)$. Then,
\begin{enumerate}

\item $ (\Upsilon_K \widehat{\bigvee} \Upsilon_L)^{\widehat{c}}=
\Upsilon_K^{\widehat{c}} \widehat{\bigwedge}
\Upsilon_L^{\widehat{c}} $
\item $ (\Upsilon_K \widehat{\bigwedge} \Upsilon_L)^{\widehat{c}} =
\Upsilon_K^{\widehat{c}} \widehat{\bigvee}
\Upsilon_L^{\widehat{c}}.$
\item  $(\Upsilon_K \widehat{\bigvee} \Upsilon_L)\widehat{\bigvee} \Upsilon_M=\Upsilon_K \widehat{\bigvee} (\Upsilon_L\widehat{\bigvee} \Upsilon_M) $
\item  $(\Upsilon_K \widehat{\bigwedge} \Upsilon_L)\widehat{\bigwedge} \Upsilon_M=\Upsilon_K \widehat{\bigwedge} (\Upsilon_L\widehat{\bigwedge} \Upsilon_M) $
\end{enumerate}
\end{prop}
\begin{pf} The proof of the Propositions 1- 4 are obvious.
\end{pf}
%---------------------------------------------------------
\begin{defn}\label{difference}
Let $\Upsilon_K, \Upsilon_L \in IVNS(U)$. Then, difference of
$\Upsilon_K$ and $\Upsilon_L$, denoted $ \Upsilon_K
\widehat{\setminus} \Upsilon_L$, is defined by
$$
\upsilon_{K \widehat{\setminus}
L}(x)=\upsilon_{K}(x)\widetilde{\setminus}\upsilon_{L}(x)\quad\textrm{for
all } x \in E.
$$
\end{defn}
%---------------------------------------------------------------
\begin{defn}\label{addition}
Let $\Upsilon_K, \Upsilon_L \in IVNS(U)$. Then, addition of
$\Upsilon_K$ and $\Upsilon_L$, denoted $ \Upsilon_K \widehat{+}
\Upsilon_L$, is defined by
$$
\upsilon_{K \widehat{+}
L}(x)=\upsilon_{K}(x)\widetilde{+}\upsilon_{L}(x)\quad\textrm{for
all } x \in E.
$$
\end{defn}
%---------------------------------------------------------
%---------------------------------------------------------------
\begin{prop}Let $\Upsilon_K, \Upsilon_L , \Upsilon_M \in IVNS(U)$. Then,
\begin{enumerate}
\item$\Upsilon_{K}(x)\widehat{+}\Upsilon_{L}(x)\widehat{=}\Upsilon_{L}(x)\widehat{+}\Upsilon_{K}(x)$
\item$(\Upsilon_{K}(x)\widehat{+}\Upsilon_{L}(x))\widehat{+}\Upsilon_{M}(x)=\Upsilon_{K}(x)\widehat{+}(\Upsilon_{L}(x)\widehat{+}\Upsilon_{M}(x))$
\end{enumerate}
\end{prop}
\begin{pf}The proofs can be easily obtained from Definition \ref{addition}.
\end{pf}
%---------------------------------------------------------------------------------
%---------------------------------------------------------------
\begin{defn}\label{multi}
Let $\Upsilon_K \in IVNS(U)$. Then, scalar multiplication of
$\Upsilon_K$, denoted $ a\widehat{\times}\Upsilon_K $, is defined by
$$
a\widehat{\times}\Upsilon_K=a\widetilde{.}\upsilon_K\quad\textrm{for
all } x \in E.
$$
\end{defn}
%---------------------------------------------------------
%---------------------------------------------------------------
\begin{prop}Let $\Upsilon_K, \Upsilon_L , \Upsilon_M \in IVNS(U)$. Then,
\begin{enumerate}
\item$\Upsilon_{K}(x)\widehat{\times}\Upsilon_{L}(x){=}\Upsilon_{L}(x)\widehat{\times}\Upsilon_{K}(x)$
\item$(\Upsilon_{K}(x)\widehat{\times}\Upsilon_{L}(x))\widehat{\times}\Upsilon_{M}(x)=\Upsilon_{K}(x)\widehat{\times}(\Upsilon_{L}(x)\widehat{\times}\Upsilon_{M}(x))$
\end{enumerate}
\end{prop}
\begin{pf}The proofs can be easily obtained from Definition
\ref{multi}.
\end{pf}
%%---------------------------------------------------------------
\begin{defn}\label{divi}
Let $\Upsilon_K \in IVNS(U)$. Then, scalar division  of
$\Upsilon_K$, denoted $ \Upsilon_K\widehat{/}a $, is defined by
$$
\Upsilon_K\widehat{/}a =\Upsilon_K\widetilde{/}a \quad\textrm{for
all } x \in E.
$$
\end{defn}
%---------------------------------------------------------
\begin{exmp}\label{ex-soft-3}
Consider the above Example \ref{ex-soft}, for $a=5$, the scalar
division  of $\Upsilon_K$, denoted $ \Upsilon_K\widehat{/}5 $, can
be represented into the following table;
$$
\begin{tabular}{|c|cc|}
  \hline
 \(U\) &  \(u_1\) &
\(u_2\)\\
    \hline \(x_1\)  &
\(\langle [0.1,0.14], [0.16,0.18], [0.04,0.1]\rangle\)& \(\langle
[0.06,0.12], [0.15,0.18], [0.04,0.16]
\rangle\) \\
    \(x_2\) &
\(\langle [0.0,0.06], [0.12,0.16], [0.16,0.18] \rangle\)& \(\langle
[0.02,0.16], [0.16,0.18], [0.15,0.25]
\rangle\) \\
    \(x_3\)  &
\(\langle [0.02,0.14], [0.08,0.1], [0.16,0.18] \rangle\)& \(\langle
[0.04,0.1], [0.1,0.14], [0.12,0.16]
\rangle\) \\
    \(x_4\)  &
\(\langle [0.04,0.08], [0.14,0.18], [0.12,0.18] \rangle\)& \(\langle
[0.15,0.18], [0.12,0.18], [0.06,0.18]
\rangle\) \\
   \(x_5\) &
\(\langle [0.0,0.04], [0.2,0.2], [0.2,0.2]\rangle\)& \(\langle
 [0.0,0.05], [0.18,0.2], [0.04,0.18]
\rangle\) \\
     \hline
\end{tabular}
$$
\begin{center}
\footnotesize{\emph{\emph{Table 9}: The tabular representation of
the $ivn$-soft set $\Upsilon_K\widehat{/}5 $}}
\end{center}
\end{exmp}
%---------------------------------------------------------------
\begin{defn}\label{truth}
Let $\Upsilon_K \in IVNS(U)$. Then, truth-Favorite  of $\Upsilon_K$,
denoted $ \widehat{\bigtriangleup} \Upsilon_K $, is defined by
$$
 \widehat{\bigtriangleup} \Upsilon_K= \widetilde{\bigtriangleup} \upsilon_K\quad\textrm{for
all } x \in E.
$$
\end{defn}
%---------------------------------------------------------
\begin{exmp}\label{ex-soft-3}
Consider the above Example \ref{ex-soft}, the truth-Favorite  of
$\Upsilon_K$, denoted $ \widehat{\bigtriangleup} \Upsilon_K $, can
be represented into the following table;
$$
\begin{tabular}{|c|cc|}
  \hline
 \(U\) &  \(u_1\) &
\(u_2\)\\
    \hline \(x_1\)  &
\(\langle [1.0,1.0], [0.0,0.0], [0.2,0.5]\rangle\)& \(\langle
[0.6,1.0], [0.0,0.0], [0.2,0.8]
\rangle\) \\
    \(x_2\) &
\(\langle [0.6,1.0], [0.0,0.0], [0.3,0.9] \rangle\)& \(\langle
[0.9,1.0], [0.0,0.0], [0.3,0.5]
\rangle\) \\
    \(x_3\)  &
\(\langle [0.5,1.0], [0.0,0.0], [0.8,0.9] \rangle\)& \(\langle
[0.7,1.0], [0.0,0.0], [0.6,0.8]
\rangle\) \\
    \(x_4\)  &
\(\langle [0.9,1.0], [0.0,0.0], [0.6,0.9] \rangle\)& \(\langle
[0.9,1.0], [0.0,0.0], [0.3,0.9]
\rangle\) \\
   \(x_5\) &
\(\langle [1.0,1.0], [0.0,0.0], [1.0,1.0]\rangle\)& \(\langle
 [0.9,1.0],[0.0,0.0], [0.2,0.9]
\rangle\) \\
     \hline
\end{tabular}
$$
\begin{center}
\footnotesize{\emph{\emph{Table 10}: The tabular representation of
the $ivn$-soft set $ \widehat{\bigtriangleup} \Upsilon_K $}}
\end{center}
\end{exmp}
%-----------------------------------------------------------
\begin{prop}Let $\Upsilon_K, \Upsilon_L \in IVNS(U)$. Then,
\begin{enumerate}
\item$  \widehat{\bigtriangleup}\widehat{\bigtriangleup}  \Upsilon_K =\widehat{\bigtriangleup}\Upsilon_K $
\item$  \widehat{\bigtriangleup} (\Upsilon_K \widehat{\cup} \Upsilon_K)\widehat{\subseteq} \widehat{\bigtriangleup} \Upsilon_K \widehat{\cup} \widehat{\bigtriangleup}\Upsilon_K$
\item$  \widehat{\bigtriangleup} (\Upsilon_K \widehat{\cap} \Upsilon_K)\widehat{\subseteq} \widehat{\bigtriangleup} \Upsilon_K \widehat{\cap} \widehat{\bigtriangleup}\Upsilon_K$
\item$  \widehat{\bigtriangleup} (\Upsilon_K \widehat{+} \Upsilon_K)=\widehat{\bigtriangleup} \Upsilon_K \widehat{+} \widehat{\bigtriangleup}\Upsilon_K$
\end{enumerate}
\end{prop}
\begin{pf}The proofs can be easily obtained from Definition \ref{union}, Definition
\ref{intersection} and Definition \ref{truth}.
\end{pf}
%-------------------------------------------------------------------------------------
\begin{defn}\label{false}
Let $\Upsilon_K \in IVNS(U)$. Then, False-Favorite  of $\Upsilon_K$,
denoted $ \widehat{\bigtriangledown} \Upsilon_K $, is defined by
$$
 \widehat{\bigtriangledown} \Upsilon_K= \widetilde{\bigtriangledown} \upsilon_K\quad\textrm{for
all } x \in E.
$$
\end{defn}
%---------------------------------------------------------
\begin{exmp}\label{ex-soft-3}
Consider the above Example \ref{ex-soft}, the False-Favorite  of
$\Upsilon_K$, denoted $ \widehat{\bigtriangledown} \Upsilon_K $, can
be represented into the following table;
$$
\begin{tabular}{|c|cc|}
  \hline
 \(U\) &  \(u_1\) &
\(u_2\)\\
    \hline \(x_1\)  &
\(\langle [0.5,0.7], [0.0,0.0],[1.0,1.0]\rangle\)& \(\langle
[0.3,0.6], [0.0,0.0], [0.5,1.0]
\rangle\) \\
    \(x_2\) &
\(\langle [0.0,0.3], [0.0,0.0], [0.9,1.0] \rangle\)& \(\langle
[0.1,0.8], [0.0,0.0], [1.0,1.0]
\rangle\) \\
    \(x_3\)  &
\(\langle [0.1,0.7], [0.0,0.0], [1.0,1.0] \rangle\)& \(\langle
[0.2,0.5], [0.0,0.0], [1.0,1.0]
\rangle\) \\
    \(x_4\)  &
\(\langle [0.2,0.4], [0.0,0.0], [1.0,1.0] \rangle\)& \(\langle
[0.3,0.9],[0.0,0.0], [0.9,1.0]
\rangle\) \\
   \(x_5\) &
\(\langle [0.0,0.2], [0.0,0.0], [1.0,1.0]\rangle\)& \(\langle
 [0.0,0.1], [0.0,0.0], [1.0,1.0]
\rangle\) \\
     \hline
\end{tabular}
$$
\begin{center}
\footnotesize{\emph{\emph{Table 11}: The tabular representation of
the $ivn$-soft set $ \widehat{\bigtriangledown} \Upsilon_K $}}
\end{center}
\end{exmp}
%-----------------------------------------------------------
\begin{prop}Let $\Upsilon_K, \Upsilon_L \in IVNS(U)$. Then,
\begin{enumerate}
\item$\widehat{\bigtriangledown}\widehat{\bigtriangledown} \Upsilon_K \widehat{=}\widehat{\bigtriangledown}\Upsilon_K $
\item$ \widehat{\bigtriangledown} (\Upsilon_K \widehat{\cup} \Upsilon_K)\widehat{\subseteq} \widehat{\bigtriangledown} \Upsilon_K \widehat{\cup} \widehat{\bigtriangledown}\Upsilon_K$
\item$  \widehat{\bigtriangledown} (\Upsilon_K \widehat{\cap} \Upsilon_K)\widehat{\subseteq} \widehat{\bigtriangledown} \Upsilon_K \widehat{\cap} \widehat{\bigtriangledown}\Upsilon_K$
\item$ \widehat{\bigtriangledown} (\Upsilon_K \widehat{+} \Upsilon_K)=\widehat{\bigtriangledown} \Upsilon_K \widehat{+} \widehat{\bigtriangledown}\Upsilon_K$
\end{enumerate}
\end{prop}
%---------------------------------------------------------------------------------
\begin{pf}The proofs can be easily obtained from Definition \ref{union}, Definition
\ref{intersection} and Definition \ref{false}.
\end{pf}
%-----------------------------------------------------------
\begin{thm}Let $P$ be the power set of all $ivn-$ soft sets defined
in the universe U. Then $(P, \widehat{\cap}, \widehat{\cup })$ is a
distributive lattice.
\end{thm}
\begin{pf}The proofs can be easily obtained by showing properties; idempotency, commutativity, associativity and distributivity
\end{pf}
%---------------------------------------------------------------------------------
\section{$ivn-$soft set based decision making}
In this section, we present an adjustable approach to $ivn-$soft set
based decision making problems by extending the approach to
interval-valued intuitionistic fuzzy soft set based decision making.
\cite{zha-13}.  Some of it is quoted from \cite{jia-11,maj-02,
wan-05, zha-13}.
%----------------------------------------------------------
\begin{defn}Let $\Upsilon_K \in IVNS(U)$. Then a \emph{relation form} of
$\Upsilon_K $ is defined by
$$
R_{\Upsilon_K}=\{({r_{\Upsilon_K}}(x,u)/(x,u)):
{r_{\Upsilon_K}}(x,u)\in IVN(U), x\in E, u\in U\}
$$
where

$r_{\Upsilon_K}:E\times U \rightarrow
IVN(U)\,\,and\,\,r_{\Upsilon_K}(x,u)=\upsilon_{K(x)}(u)$ for all $
x\in E$ and $ u\in U$.

That is, $r_{\Upsilon_K}(x,u)=\upsilon_{K(x)}(u)$ is characterized
by truth-membership function $T_K$, a indeterminacy-membership
function $I_K$ and a falsity-membership function $F_K$. For each
point $ x\in E$ and $ u\in U$; $T_K$, $I_K$ and $F_K \subseteq
[0,1]$.
\end{defn}
%---------------------------------------------------------------
\begin{exmp}\label{ex-soft-12}
Consider the above Example \ref{ex-soft}, then,
$r_{\Upsilon_K}(x,u)=\upsilon_{K(x)}(u)$ can be given as follows

$\upsilon_{K(x_1)}(u_1)= \langle [0.6,0.8], [0.8,0.9],
[0.1,0.5]\rangle$,

$\upsilon_{K(x_1)}(u_2)= \langle [0.5,0.8], [0.2,0.9], [0.1,0.7]
\rangle$,

$\upsilon_{K(x_2)}(u_1)= \langle [0.1,0.4], [0.5,0.8], [0.3,0.7]
\rangle$,

$\upsilon_{K(x_2)}(u_1)= \langle [0.1,0.9], [0.6,0.9], [0.2,0.3]
\rangle$,

$\upsilon_{K(x_3)}(u_1)= \langle [0.2,0.9], [0.1,0.5], [0.7,0.8]
\rangle$,

$\upsilon_{K(x_3)}(u_2)= \langle [0.4,0.9], [0.1,0.6], [0.5,0.7]
\rangle$,

$\upsilon_{K(x_4)}(u_1)= \langle [0.6,0.9], [0.6,0.9], [0.6,0.9]
\rangle$,

$\upsilon_{K(x_4)}(u_2)= \langle [0.5,0.9], [0.6,0.8], [0.1,0.8]
\rangle$,

$\upsilon_{K(x_5)}(u_1)= \langle [0.0,0.9], [1.0,1.0],
[1.0,1.0]\rangle$,

 $\upsilon_{K(x_5)}(u_2)= \langle
 [0.0,0.9], [0.8,1.0], [0.2,0.5]
\rangle $.
\end{exmp}
%---------------------------------------------------------------
Zhang et al.\cite{zha-13} introduced level-soft set and different
thresholds on different parameters in interval-valued intuitionistic
fuzzy soft sets. Taking inspiration these definitions we give
level-soft set and different thresholds on different parameters in
$ivn-$soft sets.
%-------------------------------------------------------------------------------------
\begin{defn}\label{level}
Let $\Upsilon_K  \in IVNS(U)$. For $\alpha, \beta, \gamma \subseteq
[0,1] $, the $(\alpha, \beta, \gamma)$-level soft set of
$\Upsilon_K$ is a crisp soft set, denoted $(\Upsilon_{K};{<\alpha,
\beta, \gamma> }) $, defined by
\end{defn}
$$
(\Upsilon_{K};{<\alpha, \beta, \gamma> })= \{(x_i,\{u_{ij}:u_{ij}\in
U, \mu(u_{ij})=1\}) : x_i\in E\}
$$
where,
$$
\mu(u_{ij})=\left\{
\begin{array}{ll}
1, & (\alpha, \beta, \gamma) \widehat{\leq }\upsilon_{K(x_i)}(u_j)\\
0, & others
\end{array}\right.
$$
for all $u_j\in U$.

Obviously, the definition is an extension of level soft sets of
interval-valued intuitionistic fuzzy soft sets \cite{zha-13}.
%---------------------------------------------------------------
\begin{rem} In Definition \ref{level}, $\alpha=(\alpha_1,\alpha_2) \subseteq
[0,1] $ can be viewed as a given least threshold on degrees of
truth-membership, $ \beta=( \beta_1, \beta_2) \subseteq [0,1] $ can
be viewed as a given greatest threshold on degrees of
indeterminacy-membership and $ \gamma=( \gamma_1, \gamma_2)
\subseteq [0,1] $ can be viewed as a given greatest threshold on
degrees of falsity-membership. If $(\alpha, \beta, \gamma)
\widehat{\leq }\upsilon_{K(x_i)}(u)$, it shows that the degree of
the truth-membership of $u$ with respect to the parameter $x$ is not
less than $\alpha$, the degree of the indeterminacy-membership of
$u$ with respect to the parameter $x$ is not more than $ \gamma$ and
the degree of the falsity-membership of $u$ with respect to the
parameter $x$ is not more than $ \beta$. In practical applications
of $inv-$soft sets, the thresholds $\alpha$, $\beta$, $\gamma$ are
pre-established by decision makers and reflect decision makers'
requirements on ''truth-membership levels'',
''indeterminacy-membership levels'' and ''falsity-membership
levels'', respectively.
\end{rem}
%---------------------------------------------------------------------------------
\begin{exmp}\label{ex-soft-15}
Consider the above Example \ref{ex-soft}.

Clearly the $( [0.3,0.4], [0.3,0.5], [0.1,0.2])$-level soft set of
$\Upsilon_K$ as follows
$$
(\Upsilon_{K};{< [0.3,0.4], [0.3,0.5], [0.1,0.2]>
})=\{(x_1,\{u_1\}),(x_4,\{u_1,u_2\})\}
$$
\end{exmp}
%---------------------------------------------------------------------------------
\begin{nnote} In some practical
applications the thresholds $\alpha, \beta, \gamma$ decision makers
need to impose different thresholds on different parameters.  To
cope with such problems, we replace a constant value the thresholds
by a function as the thresholds on truth-membership values,
indeterminacy-membership values and falsity-membership values,
respectively.
\end{nnote}
%%-------------------------------------------------------------------------------------
\begin{thm}Let $\Upsilon_K,\Upsilon_L \in IVNS(U)$. Then,
\begin{enumerate}
\item $(\Upsilon_{K};<\alpha_1, \beta_1, \gamma_1>)$ and $(\Upsilon_{K};<\alpha_2, \beta_2, \gamma_2>)$ are
$<\alpha_1, \beta_1, \gamma_1>$-level soft set and $<\alpha_2,
\beta_2, \gamma_2>$-level soft set of $\Upsilon_K$, respectively.

If $<\alpha_2, \beta_2, \gamma_2>\widehat{\leq} <\alpha_1, \beta_1,
\gamma_1> $, then we have

$(\Upsilon_{K};<\alpha_1, \beta_1, \gamma_1>)\tilde{\subseteq}
(\Upsilon_{K};<\alpha_2, \beta_2, \gamma_2>)$.

\item $(\Upsilon_{K};<\alpha, \beta, \gamma>)$ and $(\Upsilon_{L};<\alpha, \beta, \gamma>)$ are
$<\alpha, \beta, \gamma$-level soft set $\Upsilon_K$ and
$\Upsilon_L$, respectively.

If $\Upsilon_K\widehat{ \subseteq} \Upsilon_L$, then we have
$(\Upsilon_{K};<\alpha, \beta, \gamma>)\tilde{\subseteq}
(\Upsilon_{L};<\alpha, \beta, \gamma>)$.
\end{enumerate}
\end{thm}
\begin{pf} The proof of the theorems are obvious.
\end{pf}
%-------------------------------------------------------------------------------------
\begin{defn}\label{level2}
Let $\Upsilon_K  \in IVNS(U)$. Let an interval-valued neutrosophic
set $<\alpha, \beta, \gamma>_{{\Upsilon_K}}:A \rightarrow IVN(U)$ in
U which is called a threshold interval-valued neutrosophic set. The
level soft set of $\Upsilon_K$ with respect to $<\alpha, \beta,
\gamma>_{{\Upsilon_K}}$ is a crisp soft set, denoted by
$(\Upsilon_{K};<\alpha, \beta, \gamma>_{{\Upsilon_K}})$, defined by;
$$
(\Upsilon_{K};{<\alpha, \beta, \gamma> }_{{\Upsilon_K}})=
\{(x_i,\{u_{ij}:u_{ij}\in U, \mu(u_{ij})=1\}): x_i\in E\}
$$
where,
$$
\mu(u_{ij})=\left\{
\begin{array}{ll}
1, & {<\alpha, \beta, \gamma> }_{{\Upsilon_K}}(x_i) \widehat{\leq }\upsilon_{K(x_i)}(u_j)\\
0, & others
\end{array}\right.
$$
for all $u_j\in U$.
\end{defn}
Obviously, the definition is an extension of level soft sets of
interval-valued intuitionistic fuzzy soft sets \cite{zha-13}.
%-------------------------------------------------------------------------------------
\begin{rem}
In Definition \ref{level2}, $\alpha=(\alpha_1,\alpha_2) \subseteq
[0,1] $ can be viewed as a given least threshold on degrees of
truth-membership, $ \beta=( \beta_1, \beta_2) \subseteq [0,1] $ can
be viewed as a given greatest threshold on degrees of
indeterminacy-membership and $ \gamma=( \gamma_1, \gamma_2)
\subseteq [0,1] $ can be viewed as a given greatest threshold on
degrees of falsity-membership
 of $u$ with respect to the
parameter $x$.

If $ {<\alpha, \beta, \gamma> }_{{\Upsilon_K}}(x_i) \widehat{\leq
}\upsilon_{K(x_i)}(u)$ it shows that the degree of the
truth-membership of $u$ with respect to the parameter $x$ is not
less than $\alpha$, the degree of the indeterminacy-membership of
$u$ with respect to the parameter $x$ is not more than $ \gamma$ and
the degree of the falsity-membership of $u$ with respect to the
parameter $x$ is not more than $ \beta$.
\end{rem}
%-------------------------------------------------------------------------------------
\begin{defn}\label{mid}
Let $\Upsilon_K  \in IVNS(U)$. Based on $\Upsilon_K $, we can define
an interval-valued neutrosophic set $<\alpha, \beta,
\gamma>^{avg}_{{\Upsilon_K}}:A \rightarrow IVN(U)$ by
$$
<\alpha, \beta, \gamma>^{avg}_{{\Upsilon_K}}(x_i)=\sum_{u\in
U}\upsilon_{K(x_i)}(u)/{|U|}
$$
for all $x\in E$.
\end{defn}
The interval-valued neutrosophic set $<\alpha, \beta,
\gamma>^{avg}_{{\Upsilon_K}}$ is called the avg-threshold of the
$ivn-$soft set ${\Upsilon_K}$. In the following discussions, the
avg-level decision rule will mean using the avg-threshold and
considering the avg-level soft set in $ivn-$soft sets based decision
making.

Let us reconsider the $ivn-$soft set ${\Upsilon_K}$ in Example
\ref{ex-soft}. The avg-threshold $<\alpha, \beta,
\gamma>^{avg}_{{\Upsilon_K}}$ of ${\Upsilon_K}$ is an
interval-valued neutrosophic set and can be calculated as follows:
$$
<\alpha, \beta,
\gamma>^{avg}_{{\Upsilon_K}}(x_1)=\sum_{i=1}^2\upsilon_{K(x_1)}(u_i)/{|U|}=\langle
[0.55,0.8], [0.5,0.9], [0.1,0.6] \rangle
$$
$$
<\alpha, \beta,
\gamma>^{avg}_{{\Upsilon_K}}(x_2)=\sum_{i=1}^2\upsilon_{K(x_2)}(u_i)/{|U|}=\langle
[0.1,0.65], [0.55,0.85], [0.25,0.5] \rangle
$$
$$
<\alpha, \beta,
\gamma>^{avg}_{{\Upsilon_K}}(x_3)=\sum_{i=1}^2\upsilon_{K(x_3)}(u_i)/{|U|}=\langle
[0.15,0.9], [0.1,0.55], [0.6,0.75] \rangle
$$
$$
<\alpha, \beta,
\gamma>^{avg}_{{\Upsilon_K}}(x_4)=\sum_{i=1}^2\upsilon_{K(x_4)}(u_i)/{|U|}=\langle
[0.55,0.9], [0.6,0.85], [0.35,0.85] \rangle
$$
$$
<\alpha, \beta,
\gamma>^{avg}_{{\Upsilon_K}}(x_5)=\sum_{i=1}^2\upsilon_{K(x_5)}(u_i)/{|U|}=\langle
[0.0,0.9], [0.9,1.0], [0.6,0.75] \rangle
$$
Therefore, we have
$$
\begin{array}{ll}
<\alpha, \beta, \gamma>^{avg}_{{\Upsilon_K}}= &\{\langle [0.55,0.8],
[0.5,0.9], [0.1,0.6] \rangle/x_1,\langle [0.1,0.65],[0.55,0.85],\\&
[0.25,0.5] \rangle/x_2,\langle [0.15,0.9], [0.1,0.55], [0.6,0.75]
\rangle/x_3,\langle [0.55,0.9],\\& [0.6,0.85], [0.35,0.85]
\rangle/x_4, \langle [0.0,0.9], [0.9,1.0], [0.6,0.75] \rangle/x_5 \}
\end{array}
$$
%---------------------------------------------------------------------------------
\begin{exmp}\label{ex-soft-16}
Consider the above Example \ref{ex-soft}. Clearly;
$$
(\Upsilon_{K};<\alpha, \beta,
\gamma>^{avg}_{{\Upsilon_K}})=\{(x_5,\{u_2\})\}
$$
\end{exmp}
%-------------------------------------------------------------------------------------
\begin{defn}\label{Mmm}
Let $\Upsilon_K  \in IVNS(U)$. Based on $\Upsilon_K $, we can define
an interval-valued neutrosophic set $<\alpha, \beta,
\gamma>^{Mmm}_{\Upsilon_K}:A \rightarrow IVN(U)$ by
$$
\begin{array}{rl}
<\alpha, \beta, \gamma>&^{Mmm}_{\Upsilon_K}= \{<[ max_{u\in U}\{inf
T_{\upsilon_{K(x_i)}(u)}\},max_{u\in U}\{sup
T_{\upsilon_{K(x_i)}(u)}\}],\\&[  min_{u\in U}\{inf
I_{\upsilon_{K(x_i)}(u)}\},min_{u\in U}\{sup
I_{\upsilon_{K(x_i)}(u)}\}],\\&[ min_{u\in U}\{inf
F_{\upsilon_{K(x_i)}(u)}\},min_{u\in U}\{sup
F_{\upsilon_{K(x_i)}(u)}\}]
>/x_i:x_i\in E\}
\end{array}
$$
\end{defn}
%---------------------------------------------------------------
The interval-valued  neutrosophic set $<\alpha, \beta,
\gamma>^{Mmm}_{\Upsilon_K}$ is called the max-min-min-threshold of
the $ivn-$soft set $\Upsilon_K$. In what follows the Mmm-level
decision rule will mean using the max-min-min-threshold and
considering the Mmm-level soft set in $ivn-$soft sets based decision
making.
%%-------------------------------------------------------------------------------------
%\begin{defn}\label{mid}
%Let $\Upsilon_K  \in IVNS(U)$. Based on $\Upsilon_K $, we can define
%an interval-valued neutrosophic set $<\alpha, \beta,
%\gamma>^{MMM}_{\Upsilon_K}:A \rightarrow IVN(U)$ by
%$$
%\begin{array}{rl}
%<\alpha, \beta, \gamma>&^{MMM}_{\Upsilon_K}= \{<[ max_{u\in U}\{inf
%T_{\upsilon_{K(x_i)}(u)}\},max_{u\in U}\{sup
%T_{\upsilon_{K(x_i)}(u)}\}],\\&[  min_{u\in U}\{inf
%I_{\upsilon_{K(x_i)}(u)}\},min_{u\in U}\{sup
%I_{\upsilon_{K(x_i)}(u)}\}],\\&[ min_{u\in U}\{inf
%F_{\upsilon_{K(x_i)}(u)}\},min_{u\in U}\{sup
%F_{\upsilon_{K(x_i)}(u)}\}]
%>/x_i:x_i\in E\}
%\end{array}
%$$
%\end{defn}
%%---------------------------------------------------------------
%The interval-valued  neutrosophic set $<\alpha, \beta,
%\gamma>^{MMM}_{\Upsilon_K}$ is called the max-max-max-threshold of
%the $ivn-$soft set $\Upsilon_K$. In what follows the MMM-level
%decision rule will mean using the max-max-max-threshold and
%considering the MMM-level soft set in $ivn-$soft sets based decision
%making.
%-------------------------------------------------------------------------------------
\begin{defn}\label{mmm}
Let $\Upsilon_K  \in IVNS(U)$. Based on $\Upsilon_K $, we can define
an interval-valued neutrosophic set $<\alpha, \beta,
\gamma>^{mmm}_{\Upsilon_K}:A \rightarrow IVN(U)$ by
$$
\begin{array}{rl}
<\alpha, \beta, \gamma>&^{mmm}_{\Upsilon_K}= \{<[ min_{u\in U}\{inf
T_{\upsilon_{K(x_i)}(u)}\},min_{u\in U}\{sup
T_{\upsilon_{K(x_i)}(u)}\}],\\&[  min_{u\in U}\{inf
I_{\upsilon_{K(x_i)}(u)}\},min_{u\in U}\{sup
I_{\upsilon_{K(x_i)}(u)}\}],\\&[ min_{u\in U}\{inf
F_{\upsilon_{K(x_i)}(u)}\},min_{u\in U}\{sup
F_{\upsilon_{K(x_i)}(u)}\}]
>/x_i:x_i\in E\}
\end{array}
$$
\end{defn}
%---------------------------------------------------------------
The interval-valued  neutrosophic set $<\alpha, \beta,
\gamma>^{mmm}_{\Upsilon_K}$ is called the min-min-min-threshold of
the $ivn-$soft set $\Upsilon_K$. In what follows the mmm-level
decision rule will mean using the min-min-min-threshold and
considering the mmm-level soft set in $ivn-$soft sets based decision
making.

\begin{thm}Let $\Upsilon_K \in IVNS(U)$. Then, $(\Upsilon_{K};<\alpha, \beta,
\gamma>^{avg}_{{\Upsilon_K}})$, $(\Upsilon_{K};<\alpha, \beta,
\gamma>^{Mmm}_{{\Upsilon_K}}), (\Upsilon_{K};<\alpha, \beta,
\gamma>^{mmm}_{{\Upsilon_K}}), (\Upsilon_{K};<\alpha, \beta,
\gamma>^{MMM}_{{\Upsilon_K}}))$ are the avg-level soft set,
Mmm-level soft set, mmm-level soft set, MMM-level soft set of
$\Upsilon_K \in IVNS(U)$, respectively. Then,
\begin{enumerate}
    \item $(\Upsilon_{K};<\alpha, \beta,
\gamma>^{Mmm}_{{\Upsilon_K}})\tilde{\subseteq}
(\Upsilon_{K};<\alpha, \beta, \gamma>^{avg}_{{\Upsilon_K}})$
    \item $(\Upsilon_{K};<\alpha, \beta,
\gamma>^{Mmm}_{{\Upsilon_K}})\tilde{\subseteq}
(\Upsilon_{K};<\alpha, \beta, \gamma>^{mmm}_{{\Upsilon_K}})$
\end{enumerate}
\end{thm}
\begin{pf} The proof of the theorems are obvious.
\end{pf}
%%-------------------------------------------------------------------------------------
\begin{thm}Let $\Upsilon_K,\Upsilon_L \in IVNS(U)$. Then,
\begin{enumerate}
\item Let $<\alpha_1, \beta_1, \gamma_1>_{\Upsilon_K}$ and $<\alpha_2, \beta_2, \gamma_2>_{\Upsilon_K}$ be
two threshold interval-valued neutrosophic sets. Then,
$(\Upsilon_{K};<\alpha_1, \beta_1, \gamma_1>_{\Upsilon_K})$ and
$(\Upsilon_{K};<\alpha_2, \beta_2, \gamma_2>_{\Upsilon_K})$ are
$<\alpha_1, \beta_1, \gamma_1>_{\Upsilon_K}$-level soft set and
$<\alpha_2, \beta_2, \gamma_2>_{\Upsilon_K}$-level soft set of
$\Upsilon_K$, respectively.

If $<\alpha_2, \beta_2, \gamma_2>_{\Upsilon_K} \widehat{\leq}
<\alpha_1, \beta_1, \gamma_1>_{\Upsilon_K}$, then we have

$(\Upsilon_{K};<\alpha_1, \beta_1,
\gamma_1>_{\Upsilon_K})\tilde{\subseteq} (\Upsilon_{K};<\alpha_2,
\beta_2, \gamma_2>_{\Upsilon_K})$.

\item Let $<\alpha, \beta, \gamma>_{\Upsilon_K}$ be
a threshold interval-valued neutrosophic sets.

Then, $(\Upsilon_{K};<\alpha, \beta, \gamma>_{\Upsilon_K})$ and
$(\Upsilon_{L};<\alpha, \beta, \gamma>_{\Upsilon_K})$ are $<\alpha,
\beta, \gamma$-level soft set $\Upsilon_K$ and $\Upsilon_L$,
respectively.

If $\Upsilon_K\widehat{ \subseteq} \Upsilon_L$, then we have
$(\Upsilon_{K};<\alpha, \beta,
\gamma>_{\Upsilon_K})\tilde{\subseteq} (\Upsilon_{L};<\alpha, \beta,
\gamma>_{\Upsilon_K})$.
\end{enumerate}
\end{thm}
\begin{pf} The proof of the theorems are obvious.
\end{pf}

%%-------------------------------------------------------------------------------------
%\begin{defn}\label{mid}
%Let $\Upsilon_K  \in IVNS(U)$. Based on $\Upsilon_K $, we can define
%an interval-valued neutrosophic set $<\alpha, \beta,
%\gamma>^{mMM}_{\Upsilon_K}:A \rightarrow IVN(U)$ by
%$$
%\begin{array}{rl}
%<\alpha, \beta, \gamma>&^{mMM}_{\Upsilon_K}= \{<[ max_{u\in U}\{inf
%T_{\upsilon_{K(x_i)}(u)}\},max_{u\in U}\{sup
%T_{\upsilon_{K(x_i)}(u)}\}],\\&[  min_{u\in U}\{inf
%I_{\upsilon_{K(x_i)}(u)}\},min_{u\in U}\{sup
%I_{\upsilon_{K(x_i)}(u)}\}],\\&[ min_{u\in U}\{inf
%F_{\upsilon_{K(x_i)}(u)}\},min_{u\in U}\{sup
%F_{\upsilon_{K(x_i)}(u)}\}]
%>/x_i:x_i\in E\}
%\end{array}
%$$
%\end{defn}
%%---------------------------------------------------------------
%The interval-valued  neutrosophic set $<\alpha, \beta,
%\gamma>^{mMM}_{\Upsilon_K}$ is called the min-max-max-threshold of
%the $ivn-$soft set $\Upsilon_K$. In what follows the mMM-level
%decision rule will mean using the min-max-max-threshold and
%considering the mMM-level soft set in $ivn-$soft sets based decision
%making.
%---------------------------------------------------------------
Now, we construct an $ivn-$soft set decision making method by the
following algorithm;
%-----------------------------------------------------------

\emph{\textbf{ Algorithm:}}
\begin{enumerate}

\item  Input the $ivn-$soft set $\Upsilon_K$,

\item  Input a threshold interval-valued neutrosophic set $<\alpha, \beta,
\gamma>^{avg}_{\Upsilon_K}$ (or $<\alpha, \beta,
\gamma>^{Mmm}_{\Upsilon_K}, <\alpha, \beta,
\gamma>^{mmm}_{\Upsilon_K}$) by using avg-level decision rule (or
Mmm-level decision rule, mmm-level decision rule) for decision
making.

\item Compute avg-level soft set $(\Upsilon_{K};<\alpha, \beta,
\gamma>^{avg}_{{\Upsilon_K}})$ (or Mmm-level soft set
($(\Upsilon_{K};<\alpha, \beta, \gamma>^{Mmm}_{{\Upsilon_K}})$,
mmm-level soft set $(\Upsilon_{K};<\alpha, \beta,
\gamma>^{mmm}_{{\Upsilon_K}})))$

\item Present the level soft set $(\Upsilon_{K};<\alpha, \beta,
\gamma>^{avg}_{{\Upsilon_K}})$ (or the level soft set($
(\Upsilon_{K};<\alpha, \beta, \gamma>^{Mmm}_{{\Upsilon_K}}$, the
level soft set $(\Upsilon_{K};<\alpha, \beta,
\gamma>^{mmm}_{{\Upsilon_K}})))$ in tabular form.

\item Compute the choice value $c_i$ of $u_i$ for any $u_i\in U$,

\item The optimal decision is to select $u_k$ if $c_k=max_{u_i \in
U}c_i.$
\end{enumerate}
%-------------------------------------------------------
\begin{rem}If $k$ has more than one value then any one of $u_k$ may be chosen.

If there are too many optimal choices in Step 6, we may go back to
the second step and change the threshold (or decision rule) such
that only one optimal choice remains in the end.
\end{rem}
%---------------------------------------------------------------
\begin{rem}
The aim of designing the Algorithm is to solve $ivn-$soft sets based
decision making problem by using level soft sets. Level soft sets
construct bridges between $ivn-$soft sets and crisp soft sets. By
using level soft sets, we need not treat $ivn-$soft sets directly
but only cope with crisp soft sets derived from them after choosing
certain thresholds or decision strategies such as the mid-level or
the top–bottom-level decision rules. By the Algorithm, the choice
value of an object in a level soft set is in fact the number of fair
attributes which belong to that object on the premise that the
degree of the truth-membership of $u$ with respect to the parameter
$x$ is not less than ''truth-membership levels'', the degree of the
indeterminacy-membership of $u$ with respect to the parameter $x$ is
not more than ''indeterminacy-membership levels'' and the degree of
the falsity-membership of $u$ with respect to the parameter $x$ is
not more than ''falsity-membership levels''.
\end{rem}

%-------------------------------------------------------
\begin{exmp}
Suppose that a customer to select a house from the real agent. He
can construct a $ivn-$soft set $\Upsilon_K$ that describes the
characteristic of houses according to own requests. Assume that
$U=\{u_1,u_2,u_3,u_4,u_5,u_6\}$ is the universe contains six house
under consideration in an real agent and $E=\{x_1=cheap,
x_2=beatiful, x_3=green surroundings, x_4=costly, x_5= large\}$.
%-----------------------------------------------------------

Now, we can apply the method as follows:
\begin{enumerate}
\item  Input the $ivn-$soft set $\Upsilon_K$ as,
$$
\begin{tabular}{|c|cc|}
 \hline
 \(U\) &  \(u_1\) &
\(u_2\)\\
  \hline \(x_1\)  &
\(\langle [0.5,0.7], [0.8,0.9], [0.2,0.5]\rangle\)& \(\langle
[0.3,0.6], [0.3,0.9], [0.2,0.8]
\rangle\) \\
    \(x_2\) &
\(\langle [0.0,0.3], [0.6,0.8], [0.3,0.9] \rangle\)& \(\langle
[0.1,0.8], [0.8,0.9], [0.3,0.5]
\rangle\) \\
    \(x_3\)  &
\(\langle [0.1,0.7], [0.4,0.5], [0.8,0.9] \rangle\)& \(\langle
[0.2,0.5], [0.5,0.7], [0.6,0.8]
\rangle\) \\
    \(x_4\)  &
\(\langle [0.2,0.4], [0.7,0.9], [0.6,0.9] \rangle\)& \(\langle
[0.3,0.9], [0.6,0.9], [0.3,0.9]
\rangle\) \\
   \(x_5\) &
\(\langle [0.0,0.2], [1.0,1.0], [1.0,1.0]\rangle\)& \(\langle
 [0.0,0.1], [0.9,1.0], [0.2,0.9]
\rangle\) \\
     \hline
       \hline
 \(U\) &  \(u_3\) &
\(u_4\)\\
    \hline \(x_1\)  &
\(\langle [0.5,0.8], [0.8,0.9], [0.3,0.9]\rangle\)& \(\langle
[0.1,0.9], [0.5,0.9], [0.2,0.4]
\rangle\) \\
    \(x_2\) &
\(\langle [0.9,0.9], [0.2,0.3], [0.3,0.5] \rangle\)& \(\langle
[0.7,0.9], [0.1,0.3], [0.5,0.6]
\rangle\) \\
    \(x_3\)  &
\(\langle [0.8,0.9], [0.1,0.7], [0.6,0.8] \rangle\)& \(\langle
[0.8,0.9], [0.1,0.2], [0.5,0.6]
\rangle\) \\
    \(x_4\)  &
\(\langle [0.6,0.9], [0.6,0.9], [0.6,0.9] \rangle\)& \(\langle
[0.5,0.9], [0.6,0.8], [0.5,0.8]
\rangle\) \\
   \(x_5\) &
\(\langle [0.8,0.9], [0.0,0.4], [0.7,0.7]\rangle\)& \(\langle
 [0.7,0.9], [0.5,1.0], [0.6,0.5]
\rangle\) \\
     \hline
       \hline
 \(U\) &  \(u_5\) &
\(u_6\)\\
    \hline \(x_1\)  &
\(\langle [0.6,0.8], [0.8,0.9], [0.1,0.5]\rangle\)& \(\langle
[0.5,0.8], [0.2,0.9], [0.1,0.7]
\rangle\) \\
    \(x_2\) &
\(\langle [0.1,0.4], [0.5,0.8], [0.3,0.7] \rangle\)& \(\langle
[0.1,0.9], [0.6,0.9], [0.2,0.3]
\rangle\) \\
    \(x_3\)  &
\(\langle [0.2,0.9], [0.1,0.5], [0.7,0.8] \rangle\)& \(\langle
[0.4,0.9], [0.1,0.6], [0.5,0.7]
\rangle\) \\
    \(x_4\)  &
\(\langle [0.6,0.9], [0.6,0.9], [0.6,0.9] \rangle\)& \(\langle
[0.5,0.9], [0.6,0.8], [0.1,0.8]
\rangle\) \\
   \(x_5\) &
\(\langle [0.0,0.9], [1.0,1.0], [1.0,1.0]\rangle\)& \(\langle
 [0.0,0.9], [0.8,1.0], [0.2,0.5]
\rangle\) \\
     \hline
\end{tabular}
$$
\begin{center}
\footnotesize{\emph{\emph{Table 12}: The tabular representation of
the $ivn$-soft set $ \widehat{\bigtriangledown} \Upsilon_K $}}
\end{center}

\item  Input a threshold interval-valued neutrosophic set $<\alpha, \beta,
\gamma>^{avg}_{\Upsilon_K}$ by using avg-level decision rule for
decision making as;
$$
\begin{array}{ll}
<\alpha, \beta, \gamma>&^{avg}_{\Upsilon_K}= \{\langle [0.41,0.76],
[0.56,0.9], [0.18,0.63] \rangle/x_1,\langle
[0.31,0.7],[0.46,0.66],\\& [0.31,0.58] \rangle/x_2,\langle
[0.41,0.8], [0.21,0.53], [0.61,0.76] \rangle/x_3,\langle
[0.45,0.81],\\& [0.61,0.86], [0.45,0.86] \rangle/x_4, \langle
[0.25,0.65], [0.7,0.9], [0.61,0.76] \rangle/x_5 \}
\end{array}
$$

\item Compute avg-level soft set $(\Upsilon_{K};<\alpha, \beta,
\gamma>^{avg}_{{\Upsilon_K}})$ as;
$$
(\Upsilon_{K};<\alpha, \beta, \gamma>^{avg}_{{\Upsilon_K}})=\{(x_2,
\{u_3\}), (x_3,\{u_4\}), (x_4,\{u_6\}, (x_5,\{u_3\})\}
$$

\item Present the level soft set $(\Upsilon_{K};<\alpha, \beta,
\gamma>^{avg}_{{\Upsilon_K}})$ in tabular form as;
$$
\begin{tabular}{|c|cccccc|}
  \hline
 \(U\) &  \(u_1\) &
\(u_2\)&  \(u_3\) & \(u_4\)&  \(u_5\) & \(u_6\)\\
    \hline \(x_1\)  &
\(0\)& \(0\)& \(0\)& \(0\)& \(0\)& \(0\) \\
    \(x_2\) &
\(0\)& \(0\)& \(1\)& \(0\)& \(0\)& \(0\) \\
    \(x_3\)  &
\(0\)& \(0\)& \(0\)& \(1\)& \(0\)& \(0\)\\
\(x_4\)  &
\(0\)& \(0\)& \(0\)& \(0\)& \(0\)& \(1\)\\
  \(x_5\)  &
\(0\)& \(0\)& \(1\)& \(0\)& \(0\)& \(0\)\\
   \hline
    \end{tabular}
   $$
\begin{center}
\footnotesize{\emph{\emph{Table 13}: The tabular representation of
the soft set $F_X$}}
\end{center}

\item Compute the choice value $c_i$ of $u_i$ for any $u_i\in
U$ as;

$$c_1=c_2=c_5= \sum_{j=1}^5 u_{1j}= \sum_{j=1}^6 u_{2j}=
\sum_{j=1}^5
h_{5j}=0,$$

$$c_4= c_6= \sum_{j=1}^5 u_{4j}= \sum_{j=1}^6
h_{6j}=1$$

$$ c_3= \sum_{j=1}^5 u_{3j}=2$$

\item The optimal decision is to select $u_3$ since $c_3=max_{u_i \in
U}c_i.$
\end{enumerate}
\end{exmp}
%**********************************************
Note that this decision making method can be applied for group
decision making easily with help of the Definition  \ref{or} and
Definition \ref{and}.
%---------------------------------------------------------------------------------
\section{Conclusion}
In this paper, the notion of the interval valued neutrosophic soft
sets ($ivn-$soft sets) is defined which is a combination of an
interval valued neutrosophic sets\cite{wan-05}  and a soft
sets\cite{mol-99}. Then, we introduce some definitions and
operations of $ivn-$soft sets sets. Some properties of $ivn-$soft
sets which are connected to operations have been established.
Finally, we propose an adjustable approach by using level soft sets
and illustrate this method with some concrete examples. This novel
proposal proves to be feasible for some decision making problems
involving $ivn-$soft sets. It can be applied to problems of many
fields that contain uncertainty such as computer science, game
theory, and so on.

%**********************************************
\end{document}